\documentclass{amsart}
\usepackage{hyperref}%
\hypersetup{pdftex, colorlinks, citecolor=red, linkcolor=blue}
\usepackage{amsthm}
\usepackage{bbm}
\usepackage{relsize}
\usepackage{eucal}
\usepackage{url}
\usepackage{mathrsfs}
\usepackage[all]{xy}
\usepackage[perpage]{footmisc}
\usepackage{lettrine}
\usepackage[vflt]{floatflt}
\newtheoremstyle{statement}
{13pt}
{13pt}
{\sl}
{}
{\bfseries\scshape}
{.}
{.5em}
{}

\theoremstyle{statement}

\newtheorem{theorem}{Theorem}[section]
\newtheorem{lemma}[theorem]{Lemma}
\newtheorem{proposition}[theorem]{Proposition}

\newtheorem{cor}[theorem]{Corollary}

\newtheoremstyle{definition}
{13pt}
{13pt}
{}
{}
{\scshape}
{.}
{.5em}
{}

\theoremstyle{definition}

\newtheorem{definition}[theorem]{Definition}

\newtheorem{remark}[theorem]{Remark}

\newtheoremstyle{remarks}
{13pt}
{13pt}
{}
{}
{\scshape}
{.}
{.5em}
{}

\theoremstyle{remarks}

\newtheoremstyle{remarksb}
{13pt}
{13pt}
{}
{}
{\scshape}
{.}
{\newline}
{}

\theoremstyle{remarksb}

\newtheoremstyle{underlined}
{13pt}
{13pt}
{\sl}
{}
{\scshape}
{}
{.5em}
{}

\theoremstyle{underlined}

\newtheorem*{utheorem}{\underline{\textsc{Theorem}}}
\makeatletter 
\def\@seccntformat#1{\protect\makebox[0pt][r]{\@ifundefined{#1@cntformat}%
   {\csname the#1\endcsname.\quad}%
   {\csname #1@cntformat\endcsname}%
}}
\def\section@cntformat{\S\thesection.\ }
\def\subsection@cntformat{}

\makeatother

  \DeclareMathAlphabet{\mathpzc}{OT1}{pzc}{m}{it}

\begin{document}

\title{Generalizations of the primitive and normal basis theorems}


\author{Shahram Biglari}
\address{Fakult\"at f\"ur Mathematik, Universit\"at Bielefeld, D-33615, Bielefeld, Germany}
\curraddr{} \email{biglari@mathematik.uni-bielefeld.de}
\thanks{}


\subjclass[2000]{Primary 12F10 - Secondary 14A25}

\keywords{normal basis theorem, Zariski topology}

\date{2006}

\dedicatory{}

\begin{abstract}
{\smaller Using a Zariski topology associated to a finite field extensions, we give new proofs and generalize the primitive and normal basis theorems.}
\end{abstract}

\maketitle
\section{Introduction}\label{introduction}
\noindent The aim of this note is to show that many of the algebraic properties of field extensions can be restated in terms of existence of certain dominant morphisms between affine varieties. As applications of this point of view we derive generalizations of two basic theorems in algebra, namely the primitive and normal basis theorems.\\

\noindent In the first section after this introduction we prove some elementary results about Zariski topology on standard affine spaces and their subsets. All the results here are derived from definitions. We have stated the results in the forms we make use of them and not in their full generality. A reader unfamiliar with the basic definitions and exercises on the notion of Zariski topology is referred to [3, Chapter 1, \S 2].\\

\noindent In the third section which is devoted to a generalization of the primitive element theorem, we show that a field extension $E/F$ of infinite fields is separable if and only if there exist a field $K\supseteq E$ and an $F-$algebra embedding $E\to K^n$ with a Zariski dense image. This immediately implies that $E/F$ is separable if and only if ${\rm Tr}_{E/F}\neq 0$. A less immediate consequence is the following generalization of the primitive element theorem.
\begin{utheorem}[{\ref{thm:generators}}]
Let $E/F$ be a separable extension of finite degree of infinite fields and $S$ a finite set of non-constant polynomials over $F$. Then there exist infinitely many $a\in E$ such that
$$
E=F(h(a)),\quad \forall h\in S.
$$
\end{utheorem}
\noindent As shown in~\ref{norm=1}, this implies that for such an extension there are infinitely many $a\in E$ with $E=F[a]$ and $N_{E/F}(a)=1$. This immediately implies an old and well-known result: let $b,c\in {\mathbbm Q}$ with $b^2-4c\not\in {\mathbbm Q}^2$. Then the equation $x^2+bxy+cy^2=1$ has infinitely many solutions in ${\mathbbm Q}$.\\

\noindent In section four, we show that an extension $E/F$ of finite degree of infinite fields is Galois if and only if there exist an $F-$algebra embedding $E\to E^n$ with a Zariski dense image. From this we obtain the following theorem.
\begin{utheorem}[{\ref{nbt}}]
Let $E/F$ be a Galois extension of finite degree of infinite fields and $S$ a finite set of non-constant polynomials over $F$. Then there exist infinitely many $a\in E$ such that for each $h\in S$
\begin{enumerate}
\item $E=F(h(a))$ and
\item $\{\sigma(h(a))\ |\ \sigma\in Gal(E/F)\}$ is a basis of $E$ over $F$.
\end{enumerate}
\end{utheorem}
\noindent This is in fact a a simultaneously generalized primitive and normal basis theorem. Related to the classical version of these theorems, we will show in~\ref{pet+nbt} that the above implies that for an extension $E\neq F$ as in the theorem, there are infinitely many $a\in E$ with $E=F(a)$, $N_{E/F}(a)=1$, and $\{\sigma(a)\ |\ \sigma\in Gal(E/F)\}$ a basis of $E$ over $F$. Even this special case seems to be a new result.\\

\noindent We would like to mention that there are quite a few articles in literature treating various forms and generalizations of the primitive and normal basis theorem among them are [1] and [4] and the references there. 
\section{The Zariski topology for finite field extensions}\label{zariski-topology}
\noindent Let $K$ be a field and $n$ a positive integer. A subset $V\subseteq K^n$ is called Zariski closed if there is a set of polynomials $S\subseteq K[x_1,\ldots,x_n]$ such that $V$ is the set of points $(\lambda_1,\ldots,\lambda_n)\in K^n$ with $f(\lambda_1,\ldots,\lambda_n)=0$ for all $f\in S$. The complement of a Zariski closed subset is called Zariski open. It is easy to see that the collection of all Zariski open subsets of $K^n$ defines a topology and that any $K-$linear endomorphism of $K^n$ is continuous.

\begin{proposition}\label{closure}
Let $F\subseteq K$ be a field extension of infinite fields and $V$ a $F-$vector subspace of $K^n$. Then the Zariski closure of $V$ in $K^n$ is the $K-$vector subspace generated by $V$.  
\end{proposition}
\begin{proof}
Let $V_K$ be the $K-$vector subspace generated by $V$. For a chosen $K-$basis of $V_K$ from elements of $V$, we can find a $K-$linear automorphism $f$ of $K^n$ transferring $V_K$ to the Zariski closed subset $K^r$ for $r={\rm dim}_KV_K$. By definition $f$ is a homeomorphism. Since $f$ is also $F$ linear, we see that
$$
F^r\subseteq f(V) \subseteq K^r.
$$
Taking the closure, it is enough to show that $F^r$ is dense in $K^r$. For this, it suffices to show that if $p\in K[x_1,\ldots, x_r]$ vanishes on $F^r$, then $p=0$. Fix a point $(\lambda_1,\ldots,\lambda_{r-1})\in F^{r-1}$. The polynomial $p(\lambda_1,\ldots,\lambda_{r-1}, x_r)\in K[x_r]$ has infinitely many roots and hence identical to the zero polynomial. This implies that $p$ is a polynomial on $x_r$ with coefficients in $K[x_1,\ldots, x_{r-1}]$ all vanishing on the whole of $F^{r-1}$. An induction shows that $p=0$.
\end{proof}
\begin{cor}\label{irreducibility}
Every dense subset $Y$ of $K^n$, for an infinite field $K$, is irreducible, i.e. if $Y\subseteq C_1\cup\cdots\cup C_m$ for Zariski closed subsets $C_i$ of $K^n$, then $C_i=K^n$ for some $i$.
\end{cor}
\begin{proof}
Assume that $C_i\neq K^n$ for all $i\leq m$. Thus there exist polynomials $f_i\in K[x_1,\ldots,x_n]$ such that $f_i$ is zero at all points of $C_i$ but non-zero at some points of $K^n$. Consider $f:=f_1\cdots f_m$. By assumption $f=0$ as a function on $Y$ and hence on $K^n$. As in the last paragraph of the proof of~\ref{closure} we see that $f=0$ as a polynomial. Hence $f_i=0$ for some $i$. This contradiction proves the corollary.  
\end{proof}
\begin{cor}\label{reg:closure}
Let $E/F$ be an extension of finite degree $n$ of infinite fields and $K$ a field containing $F$. Then the closure of the image of the regular representation
{\setlength\arraycolsep {1pt}
\begin{displaymath}
\begin{array}{rcl}
\xi_{E/F}\colon E & \longrightarrow & {End}_F(E)\otimes_FK=M_n(K)\\
x & \longmapsto & (L_x\colon y\mapsto xy)\otimes 1.
\end{array}
\end{displaymath}}
is a $K-$vector subspace of dimension $n$.
\end{cor}
\begin{proof}
We first note that $L\colon E\to {End}_F(E)=M_n(F)$, $x\mapsto L_x$ is an $F-$algebra embedding. The assertion follows from this and~\ref{closure} by virtue of the fact that any set of $F-$linearly independent elements of $End_F(E)$ remains $K-$linearly independent in $M_n(K)$.
\end{proof}
\begin{lemma}\label{dense:n}
Let $K$ be an infinite field, $n$ a positive integer, and $p_i\in K[x_1,\cdots,x_i]$ polynomials with $\partial{p_i}/\partial{x_i}$ non-zero for all $i\leq n$. Then the mapping
$$
p\colon K^n\to K^n,\quad (\lambda_1,\ldots,\lambda_n)\mapsto (p_1(\lambda_1),\ldots,p_n(\lambda_1,\ldots,\lambda_n))
$$
has a Zariski dense image.
\end{lemma}
\begin{proof}
To simplify the notations we denote a typical element $(\lambda_1,\ldots,\lambda_n)\in K^n$ by $\bar{\lambda}$ and $p_i(\lambda_1,\ldots,\lambda_i)$ by $p_i(\bar{\lambda})$ for $i\leq n$. We prove by induction on $n$ that if $q\in K[T_1,\cdots,T_n]$ is a polynomial such that
$$
q(p_1(\bar{\lambda}),\ldots,p_n(\bar{\lambda}))=0\quad \forall \bar{\lambda}\in K^n,
$$
then $q=0$. The case $n=1$ follows from the fact that $p_1$ is not constant and hence the image of $p$, as a function on $K$, is an infinite subset so that $q(p_1(\lambda))=0$ would have infinitely many solutions, that is $q=0$. Assume that the statement above holds for $n-1$. Now let us prove the statement for $n$. We can write $q$ as a sum of $q_rT_n^r$ with $q_r\in K[T_1,\ldots,T_{n-1}]$ and $r=0,1,\ldots,m$. Let $\bar{\mu}:=(\lambda_1,\ldots,\lambda_{n-1})$ be fixed. The polynomial $p_n(\lambda_1,\ldots,\lambda_{n-1},T_n)$ is not constant and hence the set $\{p_n(\lambda_1,\ldots,\lambda_{n-1},\lambda_n)\ | \ \lambda_n\in K\}$ is infinite. This means that the polynomial $q(p_1(\bar{\mu}),\ldots,p_{n-1}(\bar{\mu}),T_n)$ has infinitely many zeros and hence for all $r$ $$q_r(p_1(\bar{\mu}),\ldots,p_{n-1}(\bar{\mu}))=0\quad \forall \bar{\mu}\in K^{n-1}.$$By induction we deduce that $q=0$.
\end{proof}
\begin{cor}\label{dense:1}
Let $K$ be an infinite field and $p$ a non-constant polynomial over $K$. Then the mapping $\tilde{p}\colon K^n\to K^n$ taking $(\lambda_1,\ldots,\lambda_n)$ to $(p(\lambda_1),\ldots,p(\lambda_n))$ has a Zariski dense image. 
\end{cor}
\section{A generalization of the primitive element theorem}
\noindent Recall that a field extension $E/F$ is said to be separable if every element $\lambda\in E$ is a simple root of a non-zero polynomial $f(x)\in F[x]$. What follows is a characterization of separable field extensions and a generalization of the primitive element theorem. For a classical proof and detail see [2, V, \S 7].
\begin{definition}\label{P_K(n)}
Let $K$ be an infinite field and $n$ a positive integer. We define $P_K(n)$ to be the set of elements $(\lambda_1,\ldots,\lambda_n)\in K^n$ with $\lambda_i\neq \lambda_j$ for all $i\neq j$.
\end{definition}
\begin{proposition}\label{sep:cri}
An extension $E/F$ of infinite fields of degree $n$ is separable if and only if there exist a field extension $K/F$ and an $F-$algebra embedding $$\overline{\xi}_{E/F}\colon E\to D_n(K)$$with Zariski dense image.
\end{proposition}
\begin{proof}
Assume that $E/F$ is separable and consider the $F-$algebra embedding $\xi_{E/F}$ defined in~\ref{reg:closure}. Since $E/F$ is separable, we see that for each element $x\in E$ the matrix $\xi_{E/F}(x)$ satisfies a separable polynomial and hence diagonalizable over any algebraic closure $K$ of $F$. But $E$ is a set of commutative elements of $M_n(K)$, and hence they can be diagonalized by the same transformation, say $P\in GL_n(K)$. In this way we obtain an $F-$algebra monomorphism $\bar{\xi}_{E/F}\colon E\to D_n(K)$ by sending $x\in E$ to $P^{-1}\xi_{E/F}(x)P$. The image of this embedding is by~\ref{reg:closure} dense. Now we prove the converse. Note that $P_K(n)$ is a non-empty open subset of $D_n(K)$. Since $\bar{\xi}_{E/F}(E)$ is dense in $D_n(K)$ we can find an element $a\in E$ such that $\bar{\xi}_{E/F}(a)\in P_K(n)$. Let $A:=\bar{\xi}_{E/F}(a)\in D_n(K)$ and note that the minimal polynomial of $A$ over $K$ is of degree $n$ and with only simple roots and dividing the minimal polynomial $g(x)$ of $a$ over $F$. This means that $g(x)$ has only simple roots and is of degree $n$. Therefore $E=F(a)$ with $a$ separable. Thus $E/F$ is separable.
\end{proof}
\begin{remark}\label{emb:rem}
The above proof shows that for a separable field extension $E/F$ of degree $n$ of infinite fields, there exists a non-singular matrix $P\in GL_n(K)$ such that $x\mapsto P^{-1}\xi_{E/F}(x)P$ is an embedding with required properties in~\ref{sep:cri}.
\end{remark}
\begin{lemma}\label{P_K(n)E}
Let $F\subseteq E\subseteq K$ be extensions of infinite fields with $E/F$ of degree $n$ and $\overline{\xi}_{E/F}\colon E\to D_n(K)$ an $F-$algebra embedding with Zariski dense image. Then
$$
\bar{\xi}_{E/F}(E)\cap P_K(n)=\{\bar{\xi}_{E/F}(a)\ |\ E=F(a)\}.
$$ 
\end{lemma}
\begin{proof}
The inclusion $\subseteq$ follows from the fact that the minimal polynomial over $F$ of any element $a\in E$ has to be divisible by the minimal polynomial over $K$ of $\bar{\xi}_{E/F}(a)$. The latter is of degree $n$ because $\bar{\xi}_{E/F}(a)\in P_K(n)$. Now we show the inclusion $\supseteq$. If $a\in E$ with $E=F(a)$ and $\bar{\xi}_{E/F}(a)\not\in P_K(n)$, then every element of $\bar{\xi}_{E/F}(E)$ will be in $K^n\setminus P_K(n)$. This is impossible because $\bar{\xi}_{E/F}(E)$ is dense in $K^n$.
\end{proof}
\begin{theorem}[{The primitive element theorem}]\label{thm:generators}
Let $E/F$ be a separable extension of finite degree of infinite fields and $S$ a finite set of non-constant polynomials over $F$. Then there exist infinitely many $a\in E$ such that
$$
E=F(h(a)),\quad \forall h\in S.
$$
\end{theorem}
\begin{proof}
The case $n:={\rm dim}_E(E)=1$ is trivial and hence we assume that $n>1$. Let $K$ and $\overline{\xi}_{E/F}\colon E\to D_n(K)$ be as in~\ref{sep:cri} and $P_K(n)$ as in~\ref{P_K(n)}. For each $h\in S$ consider $h$ as a function on $E$ and the associated morphism $\tilde{h}$ as in~\ref{dense:1}. By definition~\ref{P_K(n)}, the set $P_K(n)$ is open and non-empty and hence by~\ref{dense:1} $\tilde{h}^{-1}(P_K(n))$ is a non-empty subset of $D_n(K)$. By density of $\bar{\xi}_{E/F}(E)$ in the irreducible space $D_n(K)$ we have 
$$
\bar{\xi}_{E/F}(E)\cap \bigl(\bigcap_{h\in S} \tilde{h}^{-1}(P_K(n)) \bigr)\neq \O.
$$
Let $\bar{\xi}_{E/F}(a)$ be an element in this set and $h\in S$. Note that since $\bar{\xi}_{E/F}$ is an $F-$algebra homomorphism we have $\bar{\xi}_{E/F}\circ h=\tilde{h}\circ \bar{\xi}_{E/F}$. Therefore $\tilde{h}(\bar{\xi}_{E/F}(a))\in P_K(n)$ if and only if $\bar{\xi}_{E/F}(h(a))\in P_K(n)$. It follows from~\ref{P_K(n)E} that $E=F(h(a))$. Moreover by~\ref{irreducibility} the above intersection is an infinite set.
\end{proof}
\begin{remark}\label{rem:generators}
With the notations as in \ref{thm:generators} assume that $E\neq F$. It can in fact be shown that the set $$A:=\{aF^{\times}\ |\ E=F(h(a))\ \forall h\in S\}\subseteq E^{\times}/{F^{\times}}$$is infinite. This is particularly interesting when $S=\{x\}\subseteq F[x]$.
\end{remark}
\begin{cor}\label{norm=1}
Let $E/F$ be a non-trivial separable extension of infinite fields of finite degree. Then there exist infinitely many $a\in E$ such that $E=F(a)$ and ${N}_{E/F}(a)=1$.
\end{cor}
\begin{proof}
Consider an embedding $\bar{\xi}_{E/F}\colon E\to D_n(K)$ as in~\ref{emb:rem} so that $N_{E/F}(a)={\rm det}(\bar{\xi}_{E/F}(a))$. Using \ref{thm:generators} with $S=\{x^{n}\}$ where $n={\rm dim}_FE$ we find elements $b\in E^{\times}$ with $E=F(b^n)$. Now we can take $a=b^n{N}_{E/F}(b)^{-1}$. If there were only finitely many such $a$, say $a_1,\ldots,a_k$ corresponding to $b_1,\ldots,b_k$, then for $U=\widetilde{x^n}^{-1}(P_K(n))$ and $C=\{x\in D_n(K)\ |\ x^n-{\rm det}(x)=0\}$ we would have
$$
\bar{\xi}_{E/F}(E)\subseteq \bigl(D_n(K)\setminus U\bigr)\cup \bar{\xi}_{E/F}(b_1)C\cup\cdots\cup \bar{\xi}_{E/F}(b_k)C,
$$
thus contradicting~\ref{irreducibility}.
\end{proof}
\section{A generalization of the normal basis theorem}
\noindent Recall that a separable extension $E/F$ of finite degree of fields is Galois if the set of $F-$algebra automorphisms of $E$ has exactly $[E:F]$ elements. Similar to the results of the previous section we give a Zariski topological characterization of this notion and prove a generalized version of the normal basis theorem. For a proof of the classical version see [2, V, \S 10].
\begin{proposition}\label{gal:cri}
An extension $E/F$ of infinite fields of degree $n$ is Galois if and only if there exists an $F-$algebra embedding $$\overline{\xi}_{E/F}\colon E\to D_n(E)$$with Zariski dense image.
\end{proposition}
\begin{proof}
Let $\overline{\xi}_{E/F}\colon E\to D_n(E)$ be an $F-$algebra embedding with Zariski dense image. By~\ref{sep:cri} the extension $E/F$ is separable. For each $i=1,2,\ldots,n$ define $\pi_i\colon E\to E$ be sending $x$ to $\overline{\xi}_{E/F}(x)_{ii}$. It is clear that each $\pi_i$ is an $F-$algebra homomorphism whose kernel is by density of $\overline{\xi}_{E/F}(E)$ in $D_n(E)$ trivial. Hence each $\pi_i$ is an $F-$automorphism of $E$. Finally we note that again by density of $E$ in $D_n(E)$ the morphisms $\pi_i$'s are mutually distinct. This proves that $E/F$ is Galois. To prove the converse let $E/F$ be a finite galois extension. By~\ref{sep:cri}, there exist a field $K\supseteq E$ and an $F-$algebra embedding $\overline{\xi}_{E/F}\colon E\to D_n(K)$ with Zariski dense image. We prove that the image is contained in $D_n(E)$. As above we consider $F-$algebra monomorphisms $\pi_i\colon E\to K$. Since the extension $E/F$ is galois the image of any such homomorphism is contained in $E$. Now it is easy to see that the induced topology on $D_n(E)\subseteq D_n(K)$ is the Zariski topology.
\end{proof}
\begin{definition}\label{N_E(n)}
Let $E/F$ be a Galois extension of degree $n$ of infinite fields. We identify $Gal(E/F)$ with $\{1,\ldots,n\}$ and define $N_E(n)$ to be the set of elements $(\lambda_1,\ldots,\lambda_n)\in E^n$ with $\text{det}(\lambda_{\sigma\tau})\neq 0$.
\end{definition}
\begin{lemma}
$N_E(n)$ is open and non-empty.
\end{lemma}
\begin{proof}
Being open is trivial and the rest follows from the fact that for $\lambda:=(1,0,\ldots,0)$ we have $\text{det}(\lambda_{\sigma\tau})=\pm 1$.
\end{proof}
\begin{lemma}\label{N_E(n)E}
With notations as in~\ref{N_E(n)}, $\bar{\xi}_{E/F}(E)\cap N_E(n)$ is the set of elements $\bar{\xi}_{E/F}(a)$ such that $\{\sigma(a)\ |\ \sigma\in Gal(E/F)\}$ is a basis of $E$ over $F$.
\end{lemma}
\begin{proof}
It is enough to note that for $\lambda:=\bar{\xi}_{E/F}(a)$ we have $\lambda_{\sigma}=\sigma(a)$ for all $\sigma\in Gal(E/F)$.
\end{proof}
\begin{theorem}[{The normal basis theorem}]\label{nbt}
Let $E/F$ be a Galois extension of finite degree of infinite fields and $S$ a finite set of non-constant polynomials over $F$. Then there exist infinitely many $a\in E$ such that for each $h\in S$
\begin{enumerate}
\item $E=F(h(a))$ and
\item $\{\sigma(h(a))\ |\ \sigma\in Gal(E/F)\}$ is a basis of $E$ over $F$.
\end{enumerate}
\end{theorem}
\begin{proof}
The proof is similar to that of~\ref{thm:generators}; replace $P_K(n)$ by $P_E(n)\cap N_E(n)$ and obtain
$$
\bar{\xi}_{E/F}(E)\cap \Bigl(\bigcap_{h\in S} \tilde{h}^{-1}\bigl(P_E(n)\cap N_E(n)\bigr) \Bigr)\neq \O.
$$
Let $\bar{\xi}_{E/F}(a)$ be an element in this set and $h\in S$. By~\ref{P_K(n)E} and~\ref{N_E(n)E} the two conditions in the statement are satisfied for $a$.
\end{proof}
\begin{cor}\label{pet+nbt}
Let $E/F$ be a non-trivial Galois extension of infinite fields of finite degree. There exist infinitely many $a\in E$ such that
\begin{enumerate}
\item $E=F(a)$,
\item $N_{E/F}(a)=1$,
\item $\{\sigma(a)\ |\ \sigma\in Gal(E/F)\}$ is a basis of $E$ over $F$.
\end{enumerate}
\end{cor}
\begin{proof}
Take $S=\{x^n\}$ where $n=\text{dim}_FE$ and use the argument in~\ref{norm=1} with the same $C$ but with $U$ being replaced by $\widetilde{x^n}^{-1}\bigl(P_E(n)\cap N_E(n)\bigr)$.
\end{proof}
\nocite*
\bibliographystyle{plain}
{\footnotesize \bibliography{pn-arxiv}}
\end{document}